\documentclass[hidelinks,onefignum,onetabnum]{siamart190516}

\usepackage[ntheorem]{empheq}

\usepackage{lipsum}
\usepackage{amsfonts}
\usepackage{graphicx}
\usepackage{epstopdf}
\usepackage{algorithmic}
\ifpdf
  \DeclareGraphicsExtensions{.eps,.pdf,.png,.jpg}
\else
  \DeclareGraphicsExtensions{.eps}
\fi

\newsiamremark{remark}{Remark}
\newsiamremark{hypothesis}{Hypothesis}
\crefname{hypothesis}{Hypothesis}{Hypotheses}
\newsiamthm{claim}{Claim}

\headers{Inverse scattering by inhomogeneous media}{K. Li, B. Zhang, and H. Zhang}

\title{Reconstruction of inhomogeneous media by an iteration algorithm with a learned projector}%\thanks{Submitted to the editors DATE.
\author{Kai Li\thanks{Academy of Mathematics and Systems Science, Chinese Academy of Sciences,
Beijing 100190, China and School of Mathematical Sciences, University of Chinese Academy of Sciences,
Beijing 100049, China ({\tt likai98@amss.ac.cn})}
\and
Bo Zhang\thanks{LSEC and Academy of Mathematics and Systems Science, Chinese Academy of
Sciences, Beijing, 100190, China and School of Mathematical Sciences, University of Chinese
Academy of Sciences, Beijing 100049, China ({\tt b.zhang@amt.ac.cn})}
\and
Haiwen Zhang\thanks{Corresponding author. Academy of Mathematics and Systems Science, Chinese Academy of Sciences, Beijing 100190, China ({\tt zhanghaiwen@amss.ac.cn})}}

\usepackage{amsopn}

\usepackage{latexsym,amssymb,mathrsfs,amssymb}
\usepackage{subfigure}
\usepackage{enumerate}
\usepackage{cases}
\usepackage{empheq}

\newcommand{\R}{{\mathbb R}}

\newcommand{\N}{{\mathbb N}}
\newcommand{\C}{{\mathbb C}}

\newcommand{\be}{\begin{eqnarray}}
	\newcommand{\ben}{\begin{eqnarray*}}
		\newcommand{\en}{\end{eqnarray}}
	\newcommand{\enn}{\end{eqnarray*}}

\newcommand{\la}{\lambda}

\newcommand{\wi}{\widehat}
\newcommand{\wid}{\widetilde}

\makeatletter

\newcommand{\Rmnum}[1]{\expandafter\@slowromancap\romannumeral #1@}
\makeatother

\ifpdf
\hypersetup{
  pdftitle={Reconstruction of inhomogeneous media by an iteration algorithm with a learned projector},
  pdfauthor={Kai Li, Bo Zhang and Haiwen Zhang}
}
\fi

\begin{document}

\maketitle

% REQUIRED
\begin{abstract}
This paper is concerned with the inverse problem of reconstructing an inhomogeneous medium from the
acoustic far-field data at a fixed frequency in two dimensions.
This inverse problem is severely ill-posed (and also strongly nonlinear), and certain regularization
strategy is thus needed.
However, it is difficult to select an appropriate regularization strategy
which should enforce some a priori information of the unknown scatterer.
To address this issue, we plan to use a deep learning approach to learn some a priori information
of the unknown scatterer from certain ground truth data, which is then combined with a traditional
iteration method to solve the inverse problem.
Specifically, we propose a deep learning-based iterative reconstruction algorithm for the inverse problem,
based on a repeated application of a deep neural network and the iteratively
regularized Gauss-Newton method (IRGNM).
Our deep neural network (called the learned projector in this paper) mainly focuses on learning the a priori information of the {\em shape} of the unknown contrast with a normalization technique in the
training process and is trained to act like a projector which is helpful for
projecting the solution into some feasible region.
Extensive numerical experiments show that our reconstruction algorithm provides good reconstruction
results even for the high contrast case and has a satisfactory generalization ability.
\end{abstract}

% REQUIRED

% REQUIRED
%\textbf{MSC codes.} 68T07, 35R30, 65N21, 78A46

\section{Introduction}
\setcounter{equation}{0}

This paper is concerned with the inverse problem of scattering of time-harmonic acoustic waves from
inhomogeneous media in two dimensions.
This kind of problems arises in many applications, including sonar detection, remote sensing,
geophysical exploration, medical imaging and nondestructive testing (see, e.g., \cite{C15,C19}).

It is well-known that inverse medium scattering problems are strongly nonlinear and severely ill-posed.
Consequently, numerous iterative algorithms have been developed with various regularization strategies
for recovering inhomogeneous media (or contrasts of inhomogeneous media) from a knowledge of the
far-field data or scattered-field data. For example, a continuation method was proposed in \cite{B05}
to reconstruct the inhomogeneous medium from multi-frequency scattering data.
Specifically, the Born approximation was first used to compute the initial guess of the inhomogeneous
medium from the measured data with the lowest frequency, followed by the recursive application of
the Landweber method to obtain the reconstruction result from the measured data with multiple frequencies.
A preconditioning technique was introduced in \cite{L10} for the iteratively regularized Gauss-Newton
method (IRGNM) and applied to solve inverse medium scattering problems.
A regularized Newton method with a preconditioner was proposed in \cite{H01} for inverse medium
scattering problems, which has less computational cost compared with the standard regularized Newton method.
The contrast source inversion (CSI) method was introduced in \cite{V97} for inverse medium scattering problems.
The basic idea of CSI is to minimize the cost functional by alternatively updating the contrast and the
contrast source. A subspace-based optimization method was introduced in \cite{C09} which has similar
features with CSI. For a comprehensive discussion of regularization methods, see the
monographs \cite{E96,KNS08} and the references quoted therein.
%Recently, non-iterative algorithms have attracted much attention in inverse medium scattering problems,
%such as the linear sampling method \cite{CPP97}, the gap functional method \cite{C12}, the factorization
%method \cite{K02}, the singular sources method \cite{P05} and the approximate factorization method \cite{Q19}.
%Non-iterative methods do not need to solve the forward problem and thus are computationally
%fast. However, the reconstruction results of non-iterative algorithms are usually less
%accurate than those of iterative algorithms.

In recent years, deep learning and convolutional neural networks (CNNs) have been employed
to develop efficient methods for inverse medium scattering problems
(see, e.g., \cite{K19,L18,S19,WC19,WC18,LZR22,GRL21,LYZ22,Y19,Z20}).
For example, \cite{Y19} proposed a two-step enhanced deep learning approach, where the first step
is to retrieve the initial contrast from the scattered field by a CNN and the second step uses
a residual CNN to refine the initial contrast.
In \cite{K19}, a novel deep neural network (called SwitchNet) was proposed for solving
the inverse medium scattering problem under the small contrast assumption of the
inhomogeneous medium. SwitchNet is elaborately designed by analyzing the inherent
low-rank structure of the scattering problem and trained to map the scattered-field data to
the unknown scatterer.
In \cite{LZR22}, two physics-guided loss functions were proposed to improve the noise robustness
and the reconstruction accuracy of the deep learning approach for reconstructing the unknown
scatterer from electric far-field measurements.
In \cite{L18,WC19,WC18,Z20}, the initial contrast of the inhomogeneous medium is first retrieved
by non-CNN-based methods, and well-trained CNNs are then employed to refine
the reconstruction results. Note that the CNN in \cite{L18} was built based on
a conventional iteratively regularized algorithm.
\cite{S19} initially utilized a well-trained CNN to get a good approximation of the contrast source
and subsequently applied subspace optimization methods with the total variation regularization for
further refinement. \cite{GRL21,LYZ22} unrolled the iterative reconstruction schemes into deep neural
networks, which bridges the gap between model-based methods and data-driven deep learning methods.
For a recent review of deep learning-based approaches for inverse medium scattering problems,
see \cite{C20}.

On the other hand, deep learning has also been applied to develop effective methods for solving
inverse problems, such as computed tomography (CT) \cite{G18,J17}, magnetic resonance imaging (MRI) \cite{Y18},
optical diffraction tomography (ODT) \cite{Y20} and electrical impedance tomography (EIT) \cite{W19,WZC19}.
The reader is referred to \cite{M17,AMOS19} for a good survey of deep learning-based methods for
various inverse problems. It is worth noting that deep learning is used to learn some a prior information
of the unknown objects or a regularization functional from data for some linear inverse problems.
%some researchers have developed various iterative regularization methods based on
%deep learning methods, which benefit from both physical models and information inside existing data.
For example, \cite{Y18} proposed a novel deep learning approach for MRI which combines
the traditional compressive sensing method and a deep learning method.
The basic idea of \cite{Y18} is to unroll the alternating direction method of multipliers (ADMM) as
two deep learning architectures which automatically learn a regularization functional from training data
and produce promising results with low computational cost.
\cite{G18} proposed an approach for CT image reconstruction which combines the projected gradient descent
method and a deep learning method.
The method in \cite{G18} first trains the CNN to act like a projector onto the set of desirable
solutions and then produces the reconstructed results by incorporating this CNN into the projected
gradient descent method. This approach was extended to ODT in \cite{Y20} under the Born approximation
which is also a linear inverse problem.

The purpose of this paper is to propose an iterative regularization algorithm based on
a deep learning approach to recover the contrast of an inhomogeneous medium from the far-field data,
especially for the high contrast case.
% which is difficult to solve with very good reconstruction results.
%Here, by high contrast setting, we mean that the contrast of an inhomogeneous medium reaches a level so that researchers have great difficulties in obtaining a satisfactory reconstruction result.
%Both \cite{S19} and \cite{Z20} have made significant attempts under this setting.
Since this inverse problem is severely ill-posed, then an appropriate regularization strategy is needed
to enforce some a priori information of the unknown scatterer, which is difficult.
To cope with this issue, we plan to use a deep learning approach to learn some a priori information
of the unknown scatterer from certain ground truth data.
To this end, we reformulate our inverse problem as a regularized minimization problem with an unknown
regularization functional which is then reformulated as an equivalent constrained minimization problem
with an unknown feasible region depending on the unknown regularization functional.
%which is related to the inverse scattering problem.
Motivated by \cite{G18}, we then propose a deep learning-based algorithm to solve the constrained
minimization problem, which employs a well-trained deep neural network (called the learned projector
in this paper) to force the constraint to be satisfied and the IRGNM to minimize the data misfit term.
In our algorithm, the a priori information of the shape of the unknown scatterer is encoded in the
learned projector, which is learned directly from the ground truth data.
To cope with the high contrast case, the learned projector aims to learn the a priori information of
the shape of the unknown scatterer by using a normalization technique in the training process
(see Section \ref{S4_2}).
%In doing so, it is reasonable to assume that inhomogeneous medium with modest contrast or high contrast
%could share common shape a priori information.
We remark that the application of the normalization technique is believed to alleviate the burden
of training the deep neural networks and enhance the overall performance of the proposed algorithm.
It should be noted that the choice of the initial guesses for the unknown contrast has a crucial
impact on the convergence of the IRGNM algorithm (as discussed in Remark \ref{R2} below).
In our proposed algorithm, the learned projector is expected both to offer reliable initial guesses for
IRGNM and to be helpful for accelerating the proposed algorithm.
Extensive numerical experiments demonstrate that our algorithm has a satisfactory reconstruction
capacity and good generalization ability.

The rest of the paper is organized as follows.
Section \ref{S2} presents the direct and inverse medium scattering problems considered in the paper.
In Section \ref{S3}, the inverse problem is reformulated as an equivalent constrained minimization
problem with an unknown constraint depending on the unknown regularization functional,
and a projected iterative algorithm is then introduced for solving this problem, where the constraint
needs to be further determined.
To address this issue, we propose a learned projected iterative algorithm in Section \ref{S4},
which learns the a priori information of the shape of the unknown scatterer from ground truth data.
Numerical experiments are carried out in Section \ref{S5} to illustrate the effectiveness of our algorithm.
Some conclusions and remarks are given in Section \ref{S6}.

\section{Problem formulation}\label{S2}
\setcounter{equation}{0}

In this section, we introduce the direct and inverse medium scattering problems considered in this paper.
Precisely, let $B_\rho:=\{x\in\mathbb{R}^2:|x|<\rho\}$ with $\rho>0.$ Assume that the whole space
$\mathbb{R}^2$ is filled with an inhomogeneous medium characterized by the piecewise smooth refractive
index $n(x)>0$. Define $m(x):=n(x)-1$ to be the contrast of the inhomogeneous medium and assume that $\mathrm{supp}(m)\subset B_\rho$. We illuminate the inhomogeneous medium by the incident plane wave
$u^i=u^i(x,d):=e^{ikx\cdot d}$, where $k>0$ is the wave number and $d\in\mathbb{S}^1:=\{x\in\mathbb{R}^2:|x|=1\}$
denotes the incident direction. Then the scattering problem by the inhomogeneous medium is modeled
by the reduced wave equation
\begin{equation}\label{2.1}
	\triangle u(x) +k^2n(x)u(x)=0 \qquad \text{in}\;\;\mathbb{R}^2,
\end{equation}
where the total field $u:=u^i + u^s$ is the sum of the incident field $u^i$ and the scattered field
$u^s$, and $u^s$ is required to satisfy the Sommerfeld radiation condition
\be\label{2.2}
\lim_{r\to \infty}\sqrt{r}\bigg(\dfrac{\partial u^s}{\partial r}-iku^s\bigg) = 0, \qquad r = |x|.
\en
Moreover, it is known that the scattered field $ u^s $ has the asymptotic behavior \cite{C19}
\begin{equation*}
	u^s(x)=\dfrac{e^{ik|x|}}{\sqrt{|x|}}\biggl\{u^\infty(\hat{x}) + \mathcal{O}\bigg(\dfrac{1}{|x|}\bigg)\biggr\},
	\qquad |x|\to\infty,
\end{equation*}
uniformly for all directions $\hat{x}:= x/|x|\in\mathbb{S}^1$, where $u^\infty$ is the far-field
pattern of $u^s$. To illustrate the dependence on the direction $ d\in\mathbb{S}^1 $, we write
the far-field pattern, the scattered field and the total field as $ u^\infty(x,d)$, $u^s(x, d)$ and $u(x,d)$,
respectively. We refer to \cite{C19} for the well-posedness of the direct scattering problem (\ref{2.1})--(\ref{2.2}).
This paper considers the following inverse problem.

\textbf{Inverse problem (IP)}. Determine the contrast $m$ from the measured data
$u^\infty(\hat{x},d)$ for $\hat{x}$, $d\in\mathbb{S}^1$.

For the inverse problem (IP), we introduce the far-field operator
$\mathcal{F}:L^2(B_\rho)\to L^2(\mathbb{S}^1\times\mathbb{S}^1)$ mapping the contrast $m(x)$ to
its corresponding far-field pattern $u^\infty(\hat{x},d)$, that is,
\begin{equation*}
	\mathcal{F}(m) = u^\infty.
\end{equation*}
Note that this equation is nonlinear and severely ill-posed. For the uniqueness result of the
inverse problem (IP), we refer to \cite{C19}. In practical applications, only the noisy measured data
$u^{\infty,\delta}$ is available, where $\delta>0$ denotes the noise level (see Section \ref{S5_1_1}
for the choice of $u^{\infty,\delta}$). The present work consists in solving the perturbed equation
\begin{equation}\label{2.5}
	\mathcal{F}(m) \approx u^{\infty,\delta}
\end{equation}
for the unknown contrast $m$.

In this paper, Newton type methods are used in our numerical algorithm as backbone. In doing so,
we need the Fr\'{e}chet derivative of the far-field operator $ \mathcal{F} $.
\cite[Theorem 11.6]{C19} derived the Fr\'{e}chet derivative of the far-field operator in the 3D case.
By a similar argument as in \cite{C19}, it can be proved that the far-field operator $\mathcal{F}:m\mapsto u^\infty$
is Fr\'{e}chet differentiable, and the derivative is given by $\mathcal{F}'(m)(q)=v^\infty$,
where $q\in L^2(B_\rho)$ and $v^\infty:= v^\infty(\hat{x},d)$ is the far-field pattern of
the scattered field $v$ satisfying the inhomogeneous scattering problem
\begin{align}\label{2.6}
	\left\{
	\begin{aligned}
		&\triangle v + k^2nv = -k^2uq &&\textrm{in}\;\;\mathbb{R}^2,\\
		&\lim_{r\to \infty}\sqrt{r}\bigg(\dfrac{\partial v}{\partial r} - ikv\bigg) = 0, && r = |x|.
	\end{aligned}
	\right.
\end{align}
Here, $u=u(x,d)$ is the total field corresponding to the contrast $m$.
This means that one needs to compute the numerical solution of the scattering problem \eqref{2.6}
in order to numerically solve the Fr\'{e}chet derivative of $\mathcal{F}$.

For numerical reconstruction, it is necessary to discretize the contrast $m$.
Precisely, define $C_\rho:=[-\rho, \rho]\times [-\rho, \rho]\subset\mathbb{R}^2$
and discretize $C_\rho$ into uniformly distributed $(N\times N)$ pixels which are denoted as $x_{ij},$
$i,j=1,\ldots, N$. Then the contrast $m$ can be approximately represented by a piecewise constant,
which can be denoted by a discrete matrix
$\boldsymbol{m}=(\boldsymbol{m}_{ij})\in \mathbb{C}^{N\times N}$ with $\boldsymbol{m}_{ij}:=m(x_{ij})$
and is also called the contrast matrix in the rest of the paper.
Suppose the inhomogeneous medium is illuminated by $Q$ incident plane waves $u^i(x,d_q)$
with distinct incident directions $d_q$ ($q=1,\dots, Q$) uniformly distributed on $\mathbb{S}^1$,
and the far-field pattern is measured at $P$ distinct observation directions $\hat{x}_p$ ($p=1,\dots,P$)
uniformly distributed on $\mathbb{S}^1$.
The noisy far-field pattern $u^{\infty,\delta}$ can now be discretized as a measurement matrix
$\boldsymbol{u}^{\infty,\delta}:=(a^\delta_{p,q})\in\mathbb{C}^{P\times Q}$ with
$a^\delta_{p,q}:=u^{\infty,\delta}(\hat{x}_p,d_q),\;p=1,\dots,P,\;q=1,\dots,Q$.
Note that $\textrm{supp}(m)\subset B_\rho\subset C_\rho $.
Then the formula (\ref{2.5}) can be approximated as follows
\begin{equation}\label{1}
\boldsymbol{F}(\boldsymbol{m}) \approx \boldsymbol{u}^{\infty,\delta},
\end{equation}
where $\boldsymbol{F}$ denotes the discrete form of the far-field operator $\mathcal{F}$.
We denote the discrete form of $\mathcal{F}'$ by $\boldsymbol{F}'$.

\begin{remark}\label{R1}{\rm
The scattering problem (\ref{2.1})--(\ref{2.2}) can be numerically solved by applying the fast
Fourier transform to the well-known Lippmann-Schwinger equation in a disk containing the support
of contrast $m$, as suggested by Vainikko (see \cite{V00,H01}).
Here, the Lippmann-Schwinger equation has the following form:
\ben
u(x,d) = u^i(x,d)+ k^2\int_{\mathbb{R}^2}\Phi(x,y)m(y)u(y,d)dy,\quad x\in\mathbb{R}^2,
\enn
where $\Phi(x,y):=(i/4)H_0^{(1)}(k|x-y|),\;x\neq y$, denotes the fundamental solution to the
Helmholtz equation in two dimensions and $H_0^{(1)}$ denotes the Hankel function of the first
kind of order zero. In this paper, $\boldsymbol{F} $ and $\boldsymbol{F}'$ are implemented by
using this method with the disk to be $B_\rho$.}
\end{remark}

\section{Projected iterative algorithm} \label{S3}
\setcounter{equation}{0}

To solve the inverse problem (IP), one of the most used stable reconstruction approaches is
the variational regularization. In such case, the inverse problem (IP) can be reformulated as
the minimization problem
\begin{equation}\label{4}
	\boldsymbol{m}^*\in\underset{{\boldsymbol{m}\in\mathbb{C}^{N\times N}}}{\arg\min}
	\left|\left|\boldsymbol{F}(\boldsymbol{m})-\boldsymbol{u}^{\infty,\delta}\right|\right|^2
	+\lambda\mathcal{R}(\boldsymbol{m}),
\end{equation}
where $||\cdot||$ denotes the Frobenius norm of a matrix and $\mathcal{R}:\mathbb{C}^{N\times N}\to[0,+\infty)$
is the regularization term that encodes the a priori information about the exact contrast matrix
and penalizes unfeasible solutions, and $\lambda>0$ is the regularization parameter.
We hope to choose suitable $\mathcal{R}$ and $\lambda$ such that the exact contrast of the inverse problem (IP) solves the minimization problem \eqref{4}.
However, in practical
applications, it is difficult to determine both $\mathcal{R}$ and $\lambda$.
In this paper we try to learn $\mathcal{R}$ and $\lambda$ from the ground truth data.
To do this, we first reformulate the minimization problem \eqref{4} as a constrained minimization problem
with an unknown feasible region depending on $\mathcal{R}$ and $\lambda$
and then learn the unknown feasible region instead of the unknown $\mathcal{R}$ and $\lambda$
directly from the ground truth data, as explained below. We have the following lemma.

\begin{lemma}\label{le3.1}
	Let $\boldsymbol{m}^*_1$ be a solution of \eqref{4} and define
	$\mathcal{M}_\mathcal{R}:=\{\boldsymbol{m}\in\C^{N\times N}:\lambda\mathcal{R}(\boldsymbol{m})\le\eta\}$
	with $\eta=\lambda\mathcal{R}(\boldsymbol{m}^*_1)$. Then $\boldsymbol{m}^*_1$ is a solution of
	the constrained minimization problem
	\be\label{5}
	\boldsymbol{m}^*\in\underset{\boldsymbol{m}\in\mathcal{M}_\mathcal{R}}{\arg\min}
	\|\boldsymbol{F}(\boldsymbol{m})-\boldsymbol{u}^{\infty,\delta}\|^2.
	\en
	Moreover, a solution $\boldsymbol{m}^*_2$ of \eqref{5} also solves \eqref{4}.
\end{lemma}

\begin{proof}
	We first prove that $\boldsymbol{m}^*_1$ satisfies \eqref{5} with $\eta=\la\mathcal{R}(\boldsymbol{m}^*_1)$.
	This we do by contradiction. Suppose this is not true. Then there would exist $\wid{\boldsymbol{m}}^*_1$
	satisfying \eqref{5} and that
	\ben
	\la\mathcal{R}(\wid{\boldsymbol{m}}^*_1)\leq\eta=\la\mathcal{R}(\boldsymbol{m}^*_1),\qquad \|\boldsymbol{F}(\wid{\boldsymbol{m}}^*_1)-\boldsymbol{u}^{\infty,\delta}\|
	<\|\boldsymbol{F}(\boldsymbol{m}^*_1)-\boldsymbol{u}^{\infty,\delta}\|.
	\enn
	This means that
	\ben
	\|\boldsymbol{F}(\wid{\boldsymbol{m}}^*_1)-\boldsymbol{u}^{\infty,\delta}\|
	+\la\mathcal{R}(\wid{\boldsymbol{m}}^*_1)<\|\boldsymbol{F}(\boldsymbol{m}^*_1)
	-\boldsymbol{u}^{\infty,\delta}\|+ \la\mathcal{R}(\boldsymbol{m}^*_1),
	\enn
	which contradicts to the fact that $\boldsymbol{m}^*_1$ solves \eqref{4}.
	
	We now prove that $\boldsymbol{m}^*_2$ satisfies \eqref{4}. In fact, by the definition of
	$\boldsymbol{m}^*_2$ we have
	\ben
	\la\mathcal{R}(\boldsymbol{m}^*_2)\le\eta=\la\mathcal{R}(\boldsymbol{m}^*_1),\qquad
	\|\boldsymbol{F}(\boldsymbol{m}^*_2)-\boldsymbol{u}^{\infty,\delta}\|
	\le\|\boldsymbol{F}(\boldsymbol{m}^*_1)-\boldsymbol{u}^{\infty,\delta}\|.
	\enn
	Thus
	\ben
	\|\boldsymbol{F}(\boldsymbol{m}^*_2)-\boldsymbol{u}^{\infty,\delta}\|+\la\mathcal{R}(\boldsymbol{m}^*_2)
	\le\|\boldsymbol{F}(\boldsymbol{m}^*_1)-\boldsymbol{u}^{\infty,\delta}\|+\la\mathcal{R}(\boldsymbol{m}^*_1).
	\enn
	This, together with the definition of $\boldsymbol{m}^*_1$, implies that
	$\boldsymbol{m}^*_2$ solves \eqref{4}. The proof is thus complete.
\end{proof}
%\begin{remark}\label{R3} {\rm
%It is observed from Lemma \ref{le3.1} that the solution $ \boldsymbol{m}^*_1 $ of the minimization problem \eqref{4} belongs to $ \mathcal{M}_R $. Hence, in the rest of the paper, we assume that $\mathcal{R}$ and $\lambda$ are chosen appropriately such that the exact contrast matrix lie in $ \mathcal{M}_R $.
%}
%\end{remark}

By Lemma \ref{le3.1} we know that for given $\mathcal{R}$ and $\lambda$ the minimization problem \eqref{4} is
equivalent to the constrained minimization problem \eqref{5}.
Moreover, for heavily noised measured far-field data, solving \eqref{5} will yield
a solution $\widehat{\boldsymbol{m}}\in\mathcal{M}_\mathcal{R}$ such that
$\widehat{\boldsymbol{u}}^{\infty,\delta}:=\boldsymbol{F}(\widehat{\boldsymbol{m}})$
is as close to $\boldsymbol{u}^{\infty,\delta}$ as possible. Therefore, the constrained solution
$\boldsymbol{m}\in\mathcal{M}_\mathcal{R}$ of the constrained minimization problem \eqref{5}
may have the effect of denoising the noised measured far-field
$\boldsymbol{u}^{\infty,\delta}$, as demonstrated in the numerical experiments in Section \ref{S5_3}.

From Lemma \ref{le3.1}, instead of solving the minimization problem \eqref{4}
we propose a projected iterative algorithm to solve the constrained minimization problem \eqref{5}.
Specifically, given a projection operator $\mathcal{P}$ which projects
$\boldsymbol{m}\in\C^{N\times N}$ onto the feasible region $\mathcal{M}_R$,
we first employ the Landweber method and $\mathcal{P}$ to generate an initial guess of the unknown
contrast and then use the iteratively regularized Gauss-Newton method (IRGNM) and $\mathcal{P}$
alternately to improve the initial guess.
We will introduce the Landweber method and IRGNM in Sections \ref{S3_1} and \ref{S3_2}, respectively.
The projected iterative algorithm will be given in Section \ref{S3_3}.

It should be noted that the information of $\mathcal{R}$ and $\lambda$ in \eqref{5} is encoded in $\mathcal{P}$.
Intuitively, the regularization term $\mathcal{R}$ should be zero when meeting the ground truth contrasts,
and tend to be large when encountering unfeasible solutions. Correspondingly, $\mathcal{M}_\mathcal{R}$
should contain the ground truth contrasts and exclude unfeasible solutions.
As a result, an ideal $\mathcal{P}$ should keep $\boldsymbol{m}\in\mathcal{M}_\mathcal{R}$ unchanged and
project unfeasible solutions onto $\mathcal{M}_\mathcal{R}$. However, determining $\mathcal{P}$ is challenging.
In this paper, we will develop a deep learning approach to learn $\mathcal{P}$ from  ground truth data
in Section \ref{S4}.

\subsection{Landweber method}\label{S3_1}

The Landweber iteration has been extensively studied for linear and nonlinear ill-posed problems
(see, e.g., \cite{H95,H99,K11}) and is given as follows:
\begin{equation}\label{3.2}
	\boldsymbol{m}^\delta_{i+1} =\boldsymbol{m}^\delta_i
	+\mu[\boldsymbol{F}'(\boldsymbol{m}^\delta_i)]^*(\boldsymbol{u}^{\infty,\delta}
	-\boldsymbol{F}(\boldsymbol{m}^\delta_i)),
\end{equation}
where $\boldsymbol{m}^\delta_i$ and $\boldsymbol{m}^\delta_{i+1}$ are approximations of
the unknown contrast at the $i$-th and $(i+1)$-th iterations, respectively,
$[\boldsymbol{F}'(\boldsymbol{m}^\delta_i)]^*$ denotes the adjoint of
$\boldsymbol{F}'(\boldsymbol{m}^\delta_i)$ and $\mu>0$ is the stepsize.
%of the Landweber method.
Here, the superscript $\delta$ indicates the dependence on the noise level.
Let the initial guess to be $0$ and let $L$ be the total iteration number.
The Landweber method is presented in Algorithm \ref{L}.

\begin{algorithm}[htbp]

\caption{Landweber method}\label{L}

\textbf{Input:} $\boldsymbol{F}$, $\boldsymbol{u}^{\infty,\delta}$, $\mu$, $L$

\textbf{Output:} final approximate contrast for $\boldsymbol{m}$

\textbf{Initialize:} $i=0$, $\boldsymbol{m}^\delta_0 = 0$

\begin{algorithmic}[1]
\STATE{ \textbf{while} $i<L$ \textbf{do}}
\STATE{ \begin{itemize}
\item[]
$\boldsymbol{m}^\delta_{i+1}=\boldsymbol{m}^\delta_i + \mu[\boldsymbol{F}'(\boldsymbol{m}^\delta_i)]^*(\boldsymbol{u}^{\infty,\delta}
-\boldsymbol{F}(\boldsymbol{m}^\delta_i))$
\end{itemize}}
\STATE{ \begin{itemize}
\item[] $i\gets i+1$	
\end{itemize}}
\STATE{ \textbf{end while}}
\STATE{ Set the final approximate contrast to be $\boldsymbol{m}_L^\delta$.}
\end{algorithmic}
\end{algorithm}

\begin{remark}\label{re3.2} {\rm
		It is known that the Landweber method is very cheap. However, both theoretical analysis and
		numerical experiments show that the standard Landweber method has low convergence rate
		(see, e.g., \cite{H01,E96,H99}). Hence, in this paper, we only use Algorithm \ref{L}
		for generating an acceptable initial guess.
	}
\end{remark}

\subsection{Iteratively regularized Gauss-Newton method}\label{S3_2}

Iteratively regularized Gauss-Newton method (IRGNM) was first proposed by Bakushinskii \cite{B92},
which is an inexact Newton method that incorporates the initial guess as an important a priori
information for regularization.
Precisely, let $\boldsymbol{m}_0^I$ be the initial guess of the unknown contrast matrix $\boldsymbol{m}$,
set $\boldsymbol{m}^\delta_0:=\boldsymbol{m}_0^{I}$ and define $\boldsymbol{h}_i:=\boldsymbol{m}_{i+1}^{\delta}-\boldsymbol{m}_{i}^{\delta}$, $i=0,1,2,\ldots$,
where $\boldsymbol{m}_{i}^{\delta}$ and $\boldsymbol{m}_{i+1}^{\delta}$ are approximations to
the unknown contrast at the $i$-th and $(i+1)$-th iterations, respectively.
Then the update $\boldsymbol{h}_{i}$ can be computed as follows:
\begin{equation}\label{2}
	\boldsymbol{h}_{i} = \left(\alpha_i \boldsymbol{I} + [\boldsymbol{F}'(\boldsymbol{m}^\delta_{i})]^*\boldsymbol{F}'(\boldsymbol{m}^\delta_{i}) \right)^{-1}\left([\boldsymbol{F}'(\boldsymbol{m}^\delta_{i})]^*\left(\boldsymbol{u}^{\infty,\delta}
	-\boldsymbol{F}(\boldsymbol{m}^\delta_{i})\right) + \alpha_i(\boldsymbol{m}^\delta_{0}
	-\boldsymbol{m}^\delta_{i})\right),
\end{equation}
where $\boldsymbol{I}\in\mathbb{C}^{N\times N}$ is an identity matrix and $\{\alpha_i\}_{i=0}^\infty$
is a fixed sequence such that
\begin{equation}\label{3}
	\alpha_i >0, \quad\alpha_{i+1}\leq\alpha_i\leq\sigma\alpha_{i+1}, \quad \lim_{i\to\infty}\alpha_i=0
\end{equation}
for some $\sigma>1$ (suggested in \cite{H97}). For the convenience of later use, we stop the IRGNM
iteration by choosing a fixed total iteration number $R$.
For the IRGNM with an a posteriori stopping rule, we refer to \cite{H97}.
Now the IRGNM is presented in Algorithm \ref{I}.

\begin{remark}\label{R2} {\rm
		It is known from \cite[Chapter 4.5]{C19} that the IRGNM iteration has the regularization effect
		that prevents the iterations moving too far away from the initial guess $\boldsymbol{m}^{I}_{0}$.
		Moreover, Hohage obtained a convergence result of the IRGNM iteration in \cite[Theorem 2.3]{H97}
		for general (possibly nonlinear and ill-posed) inverse problem under some appropriate conditions,
		where it is assumed that the initial guess should be close enough to the ground truth
		(see \cite[formula (2.11)]{H97}). These arguments show that the choice of the initial guess
		$\boldsymbol{m}^{I}_{0} $ plays an important role in IRGNM.
	}
\end{remark}

\begin{algorithm}[htbp]

\caption{Iteratively regularized Gauss-Newton method}\label{I}

\textbf{Input:} $\boldsymbol{F}$, $\boldsymbol{u}^{\infty,\delta}$, $\boldsymbol{m}^{I}_0$,
$\{\alpha_i\}_{i = 0}^{R-1}$, $R$

\textbf{Output:} final approximate contrast for $\boldsymbol{m}$

\textbf{Initialize:} $i=0$, $\boldsymbol{m}^\delta_0 = \boldsymbol{m}^{I}_0$

\begin{algorithmic}[1]
\STATE{ \textbf{while} $i<R$ \textbf{do}}
\STATE{
\begin{itemize}
\item[]
$\boldsymbol{h}_{i}=\left(\alpha_i\boldsymbol{I} + [\boldsymbol{F}'(\boldsymbol{m}^\delta_{i})]^*\boldsymbol{F}'(\boldsymbol{m}^\delta_{i})\right)^{-1} \left\{[\boldsymbol{F}'(\boldsymbol{m}^\delta_{i})]^*\left(\boldsymbol{u}^{\infty,\delta}
-\boldsymbol{F}(\boldsymbol{m}^\delta_{i})\right)\right.$
\end{itemize}}
\begin{itemize}
\item[]
\qquad $\left.+ \alpha_i(\boldsymbol{m}^\delta_{0}
-\boldsymbol{m}^\delta_{i})\right\}$
\end{itemize}
\STATE{
\begin{itemize}
\item[]
$\boldsymbol{m}^\delta_{i+1} = \boldsymbol{m}^\delta_{i} + \boldsymbol{h}_{i}$
\end{itemize}}
\STATE{
\begin{itemize}
\item[] $i\gets i+1$
\end{itemize}}
\STATE{ \textbf{end while}}
\STATE{ Set the final approximate contrast to be $\boldsymbol{m}^\delta_{R}$.}
\end{algorithmic}
\end{algorithm}

\subsection{Description of the projected iterative algorithm}\label{S3_3}

Under the assumption that the exact contrast matrix $ \boldsymbol{m}\in \mathcal{M}_R $. We now describe the projected iterative algorithm.
Let $r_0$ be sufficiently large, take $N_0\in\mathbb{N}^+$ for stopping criterion
and choose $\{\alpha_i\}_{i=0}^{R-1}$ satisfying \eqref{3} for the IRGNM.
Then the projected iterative algorithm can be described by the following four steps.

\textbf{Step 1}. Use Algorithm \ref{L} with far-field data $\boldsymbol{u}^{\infty,\delta}$
to obtain an approximate contrast matrix $\widetilde{\boldsymbol{m}}:=\boldsymbol{m}_{L}^\delta$.
Denote by $\mathcal{L}$ the mapping from $\boldsymbol{u}^{\infty,\delta}$ to $\widetilde{\boldsymbol{m}}$.
Then set $\boldsymbol{m}_0:=\mathcal{P}(\widetilde{\boldsymbol{m}})$ to be the initial guess of
the projected iterative algorithm. Set $i\gets 0$ and go to Step 2.

\textbf{Step 2}.
If $i>N_0$, go to Step 4 and stop the algorithm. Otherwise, use Algorithm \ref{I} with
$\boldsymbol{m}_0^{I}=\boldsymbol{m}_{i}$, far-field data $\boldsymbol{u}^{\infty,\delta}$ and
$\{\alpha_i\}_{i = 0}^{R-1}$ to obtain an approximate contrast matrix $\boldsymbol{z}_{i+1}$.
Given the above far-field data used in IRGNM, let the mapping from $\boldsymbol{m}_{i}$ to
$\boldsymbol{z}_{i+1}$ be denoted by $\mathcal{I}$. Then we compute the error
$r_{i+1}:=\left\|\boldsymbol{F}(\boldsymbol{z}_{i+1})-\boldsymbol{u}^{\infty,\delta}\right\|$.
If $r_{i+1}>r_{i}$, go to Step 4 and stop the algorithm. Otherwise, go to Step 3.

\textbf{Step 3}.
Compute $\boldsymbol{m}_{i+1}:=\mathcal{P}(\boldsymbol{z}_{i+1})$.
Set $i\gets i+1$ and go to Step 2.

\textbf{Step 4}.
Set $\boldsymbol{m}_i$ to be the final approximation.

\begin{figure}[htpb]
	\centering
	\includegraphics[height=3.2cm]{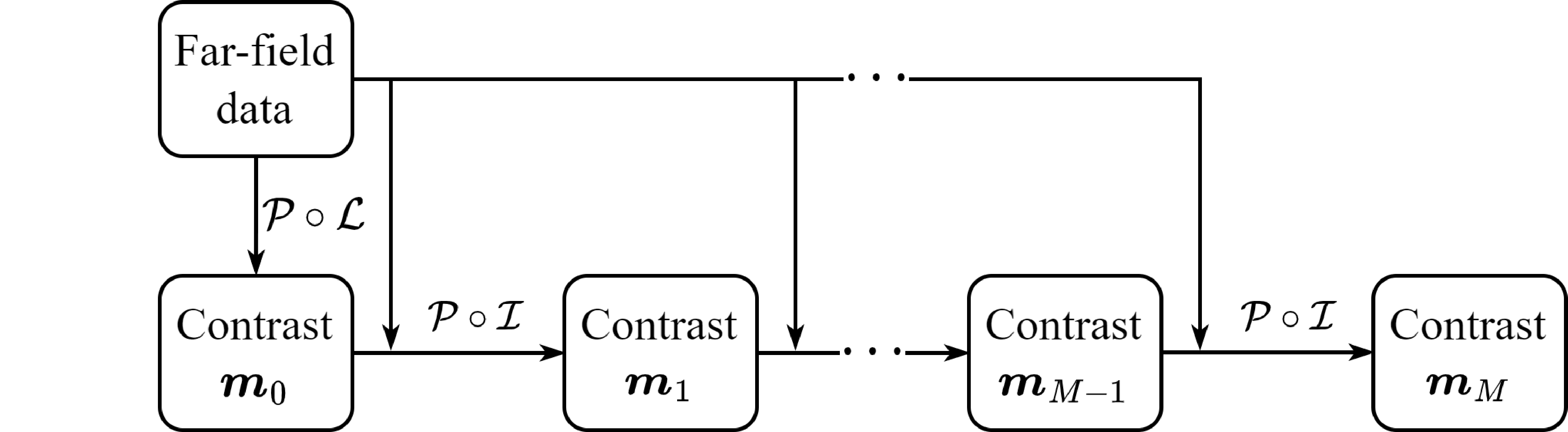}
	\caption{Diagram of the projected iterative algorithm
	}\label{F1}
\end{figure}

The above algorithm is presented in Algorithm \ref{P}. This algorithm is visualized in Figure \ref{F1},
where it is assumed that the final approximation is the output of $(M+1)$-th application of $\mathcal{P}$.
It can be seen in Algorithm \ref{P} that the IRGNM mainly contributes to the data-fitting term
$\left\|\boldsymbol{F}(\boldsymbol{m})-\boldsymbol{u}^{\infty,\delta}\right\|$ in \eqref{5}
and the projector $\mathcal{P}$ is expected to force the constraint
$\boldsymbol{m}\in\mathcal{M}_\mathcal{R}$.
Moreover, based on Remark \ref{R2} and references therein, it can be seen that the initial guess plays
an essential role in IRGNM. The projector $\mathcal{P}$ is expected to provide a good initial guess
for IRGNM and thus accelerates the iterative process of IRGNM.
Hence, it is reasonable to use Algorithm \ref{P} to solve the constrained minimization problem \eqref{5}.
%Note that we have not yet given the explicit form of $\mathcal{P}$.
%As mentioned before,
However, $\mathcal{P}$ is unknown, and we will train a deep neural network that learns the a priori
information directly from the ground truth data to provide a reasonable choice for $\mathcal{P}$
in the next section, and subsequently propose a projected iterative algorithm with a learned projector
to solve the problem (IP).

%\hphantom
\begin{algorithm}[htbp]

\caption{Projected iterative algorithm}\label{P}
\textbf{Input: } $\boldsymbol{F},\boldsymbol{u}^{\infty,\delta}$, $\mu$, $L$, $\{\alpha_i\}_{i=0}^{R-1}$,
$R$, $\mathcal{P}$, $N_0$

\textbf{Output:} final approximate contrast for $\boldsymbol{m}$

\textbf{Initialize:} $i=0$, $r_0=10^4$
\begin{algorithmic}[1]

\STATE{
Use Algorithm \ref{L} with $\boldsymbol{u}^{\infty,\delta}$ to obtain $\widetilde{\boldsymbol{m}}$
and set $\boldsymbol{m}_0:=\mathcal{P}(\widetilde{\boldsymbol{m}})$ to be the initial guess.}

\STATE{ \textbf{while $i\leq N_0$ do}}

\STATE{
\begin{itemize}
\item[]
Use Algorithm \ref{I} with $\boldsymbol{u}^{\infty,\delta}$ and $\boldsymbol{m}_{i}$ to obtain
an approximate contrast matrix $\boldsymbol{z}_{i+1}$, then compute
$r_{i+1}:=\left\|\boldsymbol{F}(\boldsymbol{z}_{i+1})-\boldsymbol{u}^{\infty,\delta}\right\|$.
\end{itemize}}

\STATE{
\begin{itemize}
\item[]  \textbf{if $r_{i+1} > r_{i}\;$ then}
\end{itemize}}

\STATE{
\begin{itemize}
\item[]
\begin{itemize}
\item[]
Set $ \boldsymbol{m}_i $ to be the final approximation and stop the algorithm.
\end{itemize}
\end{itemize}}

\STATE{
\begin{itemize}
\item[]  \textbf{end if}
\end{itemize}}

\STATE{
\begin{itemize}
\item[]
Set $\boldsymbol{m}_{i+1}:=\mathcal{P}(\boldsymbol{z}_{i+1})$
\end{itemize}}

\STATE{
\begin{itemize}
\item[] $i\gets i+1$
\end{itemize}}

\STATE{ \textbf{end while}}

\STATE{ Set $\boldsymbol{m}_i$ to be the final approximation.}
\end{algorithmic}
\end{algorithm}

\section{Learned projected iterative algorithm}\label{S4}
\setcounter{equation}{0}

%As discussed in Section \ref{S3}, determining both $\mathcal{R}$ and $\lambda$ for variational
%regularization \eqref{4} can be difficult in practical applications.
%Intuitively, the choice of regularization should be selected based on the specific requirements of each
%application scenario, which inspires us to learn a priori information from the ground truth data.
%To this aim, we propose a projected iterative algorithm (i.e., Algorithm \ref{P}) to solve the
%constrained minimization problem \eqref{5}, which is related to the variational regularization.
%In Algorithm \ref{P}, the a priori information is encoded in the operator $\mathcal{P}$.
In this section, we propose a learned projector $\mathcal{P}_\Theta$ to learn certain a priori information
from the ground truth data which can provide a reasonable choice for $\mathcal{P}$.
Here, $\mathcal{P}_\Theta$ is a deep convolutional neural network with parameters $\Theta$,
which are determined during the training process.
In order to deal with the high contrast case, $\mathcal{P}_\Theta$ focuses on learning the a priori
information of the shape of the unknown contrast.
%which is not related to the value of the contrast
%As discussed before, it is believed that the inhomogeneous medium with the modest contrast or high
%contrast could share common shape a priori information.
We further propose a learned projected iterative algorithm that incorporates the learned projector
$\mathcal{P}_\Theta$ into Algorithm \ref{P}.
The architecture and the training strategy of $\mathcal{P}_\Theta$ will be given in Sections \ref{S4_1}
and \ref{S4_2}, respectively, and the description of the learned projected iterative algorithm is
presented in Section \ref{S4_3}.

\subsection{Network architecture of the learned projector}\label{S4_1}

\begin{figure}[htpb]
	\centering
	\includegraphics[width=13cm]{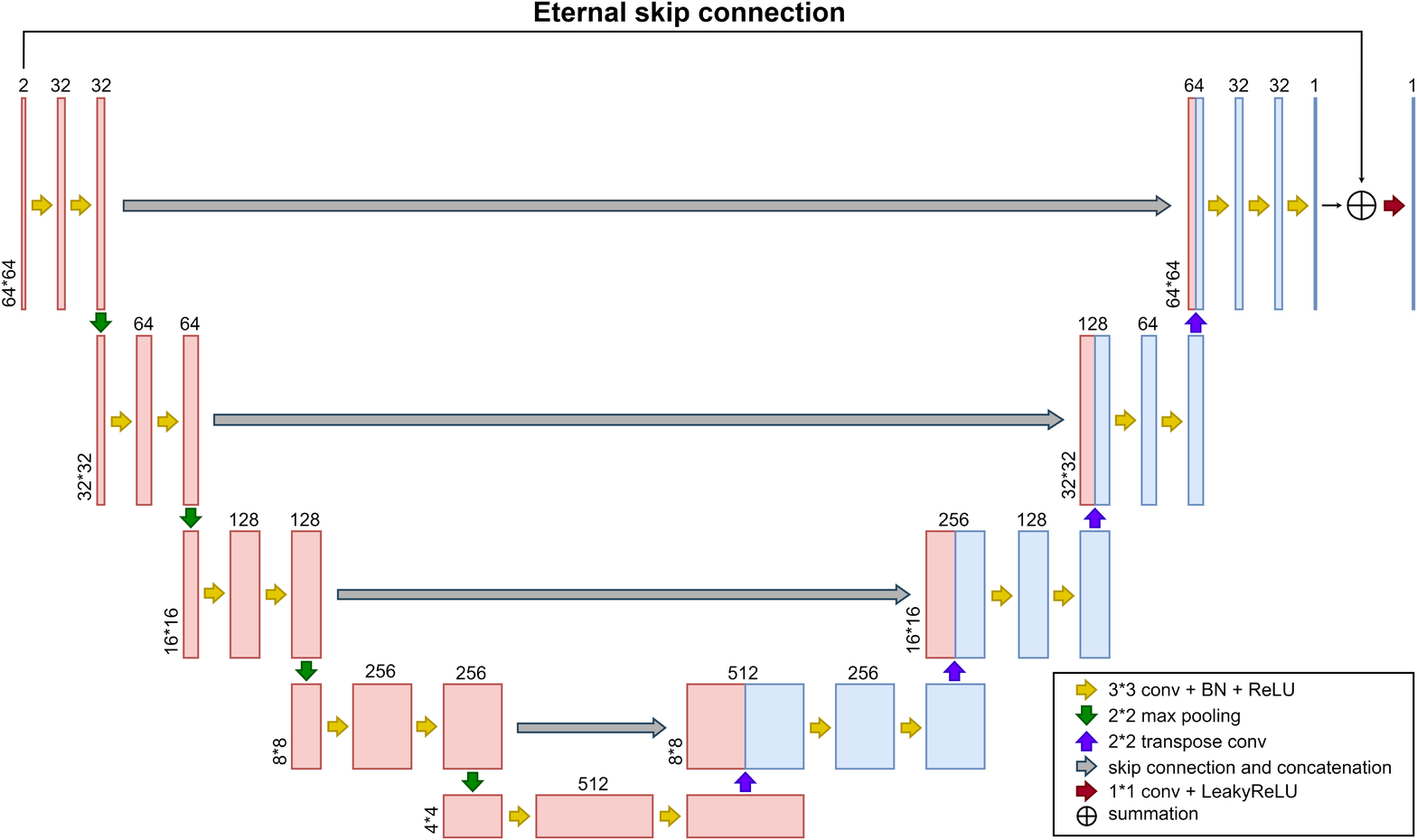}
	\caption{The architecture of $\mathcal{P}_\Theta$. Each red and blue item represents a volume
		(also called multichannel feature map). The number of channels is shown at the top of the volume,
		and the length and width are provided at the lower-left edge of the volume.
		The arrows denote different operations, which are explained at the lower-right corner of the figure.
	}\label{F2}
\end{figure}

The projector $\mathcal{P}_\Theta$ is parameterized by a convolutional neural network called U-Net,
which has a U-shaped structure. The original version of U-net was first proposed in \cite{R15}
for biomedical image segmentation. In this paper, we adopt a modified version of U-net
(see Figure \ref{F2}), which is similar to the one used in \cite{Y20}.
To be more specific, the input and the output of $\mathcal{P}_\Theta$ are the volumes with the sizes
$(N\times N\times 2)$ and $(N\times N\times 1)$, respectively,
where $N$ is given as in Section \ref{S2}.
For our proposed algorithm in Section \ref{S3}, we choose $N=64$.
%and the first channel and the second channel of the input of $\mathcal{P}_\Theta$
%will be fed with the real part and the imaginary part of the approximated contrast matrix, respectively.
%Since it is assumed in Section \ref{section2} that the exact contrast is real-valued, we set the output
%of $ \mathcal{P}_\Theta $ to be a volume with size $ (N_1\times N_1 \times 1) $, which represents the
%updated approximated contrast matrix.
As shown in Figure \ref{F2}, each red and blue item represents a
volume (also called multichannel feature map \cite{Y20}), the number of channels is shown at the top
of the volume, and the length and width are provided at the lower-left edge of the volume. The left part
and the right part of $ \mathcal{P}_\Theta $ are the contracting path and the expansive path, respectively.
For each convolutional layer of these two paths, we employ a $(3\times 3)$
convolution with zero-padding and $(1\times1)$ convolution stride, batch normalization (BN),
and rectified linear unit (ReLU) (see yellow right arrow in Figure \ref{F2}). For each down-sampling
layer in the contracting path, we apply a $ (2 \times 2) $ max pooling layer (see green downward arrow
in Figure \ref{F2}). For each up-sampling layer in the expansive path, we use a $(2\times 2)$
transposed convolution (see purple upward arrow in Figure \ref{F2}). Each up-sampled output in
the expansive path is concatenated with the corresponding multichannel feature map from the
contracting path (see gray right arrow in Figure \ref{F2}).
Moreover, the first channel of the input is added to the output of the penultimate layer (see the external
skip connection in Figure \ref{F2}).
At last, we add a Leaky rectified linear unit (LeakyReLU) \cite{M13} behind the $(1\times1)$
convolution with $(1\times1)$ convolution stride to obtain the final output (see red right
arrow in Figure \ref{F2}). For more details of U-net, see \cite{R15,J17}.

\subsection{Training strategy of the learned projector}\label{S4_2}

As discussed in Section \ref{S3}, the exact contrast matrices are assumed to lie in
$\mathcal{M}_\mathcal{R}\subset\R^{N\times N}$.
Define $\mathcal{M}_0:=\{\mathcal{N}(\boldsymbol{m}):\boldsymbol{m} = (\boldsymbol{m}_{ij})\in\mathcal{M}_\mathcal{R}\; \textrm{with}\; \boldsymbol{m}_{ij}\neq 0,\; 1\leq i,j\leq N\}$,
where $\mathcal{N}$ is a normalization operator defined by
$\mathcal{N}(\boldsymbol{f}):=\boldsymbol{f}/\|\boldsymbol{f}\|_{\max}$ for any
$\boldsymbol{f} := (\boldsymbol{f}_{ij})\in\C^{N\times N}$ with the norm $\|\boldsymbol{f}\|_{\max} := \max_{1\leq i,j\leq N}|\boldsymbol{f}_{ij}| $.
It is believed that the application of $\mathcal{N}$ over a contrast matrix $\boldsymbol{m}$ can characterize
the shape of $\boldsymbol{m}$, and thus $\mathcal{M}_0$ contains all possible shapes of the exact contrast
matrices that we are interested in (including the high contrast case).
We hope to obtain a suitable $\mathcal{P}_\Theta$ by the training process so that such a learned projector
$\mathcal{P}_\Theta$ could learn the a priori information of the shape of the unknown contrasts we are
interested in and force the normalization of approximate contrast matrices lie in $\mathcal{M}_0$.
It should be noted that $\mathcal{M}_\mathcal{R}$ is highly correlated with $\mathcal{M}_0$.
In fact, if $\boldsymbol{m}\in\mathcal{M}_\mathcal{R}$ implies that
$c\cdot\boldsymbol{m}\in\mathcal{M}_\mathcal{R}$ for any $ c\in\R$,
then $\mathcal{N}(\boldsymbol{m})\in\mathcal{M}_0$ implies that $\boldsymbol{m}\in\mathcal{M}_\mathcal{R}$.

We now describe the training strategy of the learned projector $\mathcal{P}_\Theta$
by using the normalization technique.
In the training stage, we generate a sample set of the exact contrast matrices
$\{\boldsymbol{m}^{(i)}\}_{i=1}^T$ with $T\in\N^+$ and $\boldsymbol{m}^{(i)}\in\R^{N\times N}$
such that $\|\boldsymbol{m}^{(i)}\|_{\max},\;\;i= 1,\ldots,T$, lie uniformly in the
interval $[a,b]$ with $0<a<b$ (see Section \ref{S5_2} for the choice of the interval $[a,b]$).
Note that the exact contrast matrices $\boldsymbol{m}^{(i)}$ ($i=1,\ldots,T$) are all chosen to be real
since it is assumed in Section \ref{S2} that the exact contrast is real-valued.
In what follows, let the input/output pair $(x,y)$ represent the sample of any labeled dataset,
where $x$ and $y$ denote the input  and the output, respectively.
For the input/output pair $(x,y)$ used later, $x$ and $y$ will be chosen to be a complex matrix
and a real matrix, respectively. Moreover, the real part and the imaginary part of the matrix $x$
will be put into the first channel and the second channel of the input of
$\mathcal{P}_\Theta$, respectively, and the matrix $y$ will be used as the output of
$\mathcal{P}_\Theta$. Now the training process can be
divided into the following two parts.

\textbf{Part I}. We use the approximate contrast matrices generated by the Landweber method to
train $\mathcal{P}_\Theta$.
Precisely, for each exact contrast matrix $\boldsymbol{m}^{(i)}$ ($i=1,\ldots,T$), we generate
the corresponding output of the mapping $\mathcal{L}$ (see Step 1 in Section \ref{S3_3}),
which is denoted as $\boldsymbol{m}_1^{(i)}$.
Following the training strategy in \cite{G18} and applying the normalization operator $\mathcal{N}$,
we divide this part into the following three steps.

\textbf{Step 1}. Train $ \mathcal{P}_\Theta $ on the dataset $S_1$ with Xavier initialization \cite{GB10},
where
\ben
S_1:=\{(\boldsymbol{l}_1^{(i)},\mathcal{N}(\boldsymbol{m}^{(i)}))\}_{i = 1}^T
\enn
with $\boldsymbol{l}_1^{(i)}:=\mathcal{N}(\boldsymbol{m}_1^{(i)}),\;i=1,\dots,T$.
Then we obtain $ \mathcal{P}_{\Theta_1}$ after $t_1$ epochs. After this step, it is hoped that
$\mathcal{P}_\Theta$ could learn to force the normalization of an approximate contrast matrix lie
in $\mathcal{M}_0$ (intuitively, refine the shape of the approximate contrast matrix).

\textbf{Step 2}. Train $\mathcal{P}_{\Theta}$ on the dataset $S_1\cup S_2$ with $\Theta_1$ being
the initial values of $\Theta$, where
\ben
S_2:=\{(\boldsymbol{l}_2^{(i)},\mathcal{N}(\boldsymbol{m}^{(i)}))\}_{i = 1}^T
\enn
with $\boldsymbol{l}_2^{(i)}:=\mathcal{P}_{\Theta_1}(\mathcal{N}(\boldsymbol{m}_1^{(i)})),\; i=1,\dots,T$.
Then we obtain $\mathcal{P}_{\Theta_2}$ after $t_2$ epochs.
We hope that the dataset $S_2$ can be helpful for training $\mathcal{P}_\Theta$ to mimic the important
property of any projector $\mathcal{P}$, that is,
$\mathcal{P}\circ\mathcal{P}=\mathcal{P}$, and thus increasing the stability of $\mathcal{P}_\Theta$.

\textbf{Step 3}. Train $\mathcal{P}_{\Theta}$ on the dataset $S_1\cup S_2\cup S_3$ with $\Theta_2$ being
the initial values of $\Theta$, where
\begin{equation*}
	S_3:=\{(\boldsymbol{l}_3^{(i)},\mathcal{N}(\boldsymbol{m}^{(i)}))\}_{i = 1}^T
\end{equation*}
with $\boldsymbol{l}_3^{(i)}:=\mathcal{N}(\boldsymbol{m}^{(i)}),\;i =1,\dots,T$.
Then we obtain $\mathcal{P}_{\Theta_3}$ after $t_3$ epochs. We hope that the dataset $S_3$
used in this step can be helpful for training $\mathcal{P}_\Theta$ to project an element in
$\mathcal{M}_0$ into itself.

In Steps 1, 2 and 3, we use the following error function
\ben
\mathcal{E}_M(\Theta):=\sum_{j=1}^{M}\sum_{i=1}^T\left\|\mathcal{P}_\Theta(\boldsymbol{l}_j^{(i)})
- \mathcal{N}(\boldsymbol{m}^{(i)})\right\|^2
\enn
with $M = 1,2,3$, respectively.

\textbf{Part II}. With the aid of the normalization operator $\mathcal{N}$, we use the approximate
contrast matrices generated by the IRGNM to train $\mathcal{P}_\Theta$.
For each exact contrast matrix $\boldsymbol{m}^{(i)}$ ($i=1,\ldots,T$), we generate the corresponding
output of the mapping $\mathcal{I}$ (see Step 3 in Section \ref{S3_3}) with the initial guess
$\|\boldsymbol{m}_1^{(i)}\|_{\max}\mathcal{P}_{\Theta_3}(\mathcal{N}(\boldsymbol{m}_1^{(i)}))$,
which is denoted as $\boldsymbol{r}^{(i)}$. Note that, according to the training process
in Part I, the normalization of $\|\boldsymbol{m}_1^{(i)}\|_{\max}\mathcal{P}_{\Theta_3}(\mathcal{N}(\boldsymbol{m}_1^{(i)}))$ is
expected to lie in $\mathcal{M}_0$ and to be closer to the exact contrast matrix ${\boldsymbol{m}}^{(i)}$,
compared with $\boldsymbol{m}_1^{(i)}$. Next, we train
$\mathcal{P}_{\Theta}$ on the dataset $\bigcup^4_{j=1}S_j$ with $\Theta_3$ being the initial values
of $\Theta$, where
\ben
S_4:=\{(\boldsymbol{l}_4^{(i)}, \mathcal{N}(\boldsymbol{m}^{(i)}))\}_{i = 1}^T
\enn
with $\boldsymbol{l}_4^{(i)}:=\mathcal{N}(\boldsymbol{r}^{(i)}),\;i=1,\dots,T$.
The error function in this part is given as follows:
\ben
\mathcal{E}(\Theta):=\sum_{j=1}^{4}\sum_{i=1}^T\left\|\mathcal{P}_\Theta(\boldsymbol{l}_j^{(i)})
-\mathcal{N}(\boldsymbol{m}^{(i)})\right\|^2.
\enn
Finally, we obtain the learned projector $\mathcal{P}_{\widehat{\Theta}}$ after $t_4$ epochs, which
%The learned projector $\mathcal{P}_{\widehat{\Theta}}$
is used in our proposed algorithm in the next section.
By using the dataset $S_4$ in Part II, we hope that $\mathcal{P}_\Theta$ can force the normalization
of the approximate contrast matrix generated by IRGNM lie in $\mathcal{M}_0$.
In summary, all the training datasets $S_i,\;i=1,\ldots,4,$ used in this section are expected
to make $\mathcal{P}_\Theta$ behave like a true projector onto $\mathcal{M}_0$ and to be suitable
for the projected iterative algorithm proposed in Section \ref{S3_3}.

\subsection{Description of learned projected iterative algorithm}\label{S4_3}

For any matrix $\boldsymbol{m}\in\C^{N\times N}$ define $\mathcal{A}_{\wi{\Theta}}(\boldsymbol{m})
:=\|\boldsymbol{m}\|_{\max}\mathcal{P}_{\wi{\Theta}}(\mathcal{N}(\boldsymbol{m})).$
By the property of the learned projector $\mathcal{P}_{\wi{\Theta}}$, and
according to the training strategy in Section \ref{S4_2}, the normalization of
$\mathcal{A}_{\wi{\Theta}}(\boldsymbol{m})$ is expected to lie in $\mathcal{M}_0$
and also to be closer to the ground truth contrast matrix, compared with $\boldsymbol{m}$ itself.
Moreover, if the feasible set $\mathcal{M}_\mathcal{R}$ in \eqref{5} satisfies the condition
in the first paragraph of Section \ref{S4_2}, which coincides with the high contrast case,
then $\mathcal{A}_{\widehat{\Theta}}(\boldsymbol{m})\in\mathcal{M}_\mathcal{R}$.
Therefore, the proposed reconstruction algorithm is given by Algorithm \ref{P} with $\mathcal{P}$ replaced
by $\mathcal{A}_{\widehat{\Theta}}$. For simplicity, we also called the proposed algorithm as
\textit{Learned Projected Algorithm} in the rest of the paper. See Section \ref{S5} for the
performance of this algorithm.
To show the advantage offered by the training dataset $S_4$ for the proposed algorithm,
we define
$\mathcal{A}_{\Theta_3}(\boldsymbol{m}):=\|\boldsymbol{m}\|_{\max}\mathcal{P}_{\Theta_3}
(\mathcal{N}(\boldsymbol{m}))$ for $\boldsymbol{m}\in\C^{N\times N}$.
Here, $S_4$ and $\mathcal{P}_{\Theta_3}$ are given as in Section \ref{S4_2}.
In the numerical experiments, we also consider a simplified version of \textit{Learned Projected Algorithm},
that is, Algorithm \ref{P} with $\mathcal{P}$ replaced by $\mathcal{A}_{\Theta_3}$.
This simplified algorithm is called \textit{Simplified Learned Projected Algorithm}
in the rest of the paper.

\section{Numerical experiments}\label{S5}
\setcounter{equation}{0}

In this section, we present numerical experiments to illustrate the effectiveness of the learned projected
iterative algorithm (i.e., \textit{Learned Projected Algorithm}) for the problem (IP).
The experimental setup is given in Subsection \ref{S5_1} for numerical experiments.
The performances of \textit{Learned Projected Algorithm} and \textit{Simplified Learned Projected Algorithm}
are shown in Subsection \ref{S5_2}.
To show the robustness of our algorithms with respect to noise, our algorithms are tested with different noise
levels in Subsection \ref{S5_3}.

\subsection{Experimental setup}\label{S5_1}

The training process is performed on COLAB (Tesla P100 GPU, Linux operating system) and is implemented
in PyTorch, while the computation of the direct scattering problem, the Landweber method and IRGNM are
implemented in Python 3.7 on a desktop computer (Intel Core i7-10700 CPU (2.90 GHz), 32 GB of RAM,
Ubuntu 20.04 LTS).

\subsubsection{Simulation setup for the scattering model}\label{S5_1_1}

As mentioned in Section \ref{S2}, the support of the unknown contrast is assumed to lie in
a disk $B_\rho\subset C_\rho$ with $\rho>0$. Without loss of generality, we choose $\rho=1$.
The number of incident directions is set to be $Q=16$ and the number of measured directions is $P=32$.
To generate the synthetic far-field data, we use the method discussed in Remark \ref{R1} with $N=256$.
The noisy far-field data $\boldsymbol{u}^{\infty,\delta}(\hat{x}_p,d_q)$,
$p=1,\ldots,P,\; q =1,\ldots,Q$, are given as
\ben
\boldsymbol{u}^{\infty,\delta}(\hat{x}_p,d_q)=\boldsymbol{u}^{\infty}(\hat{x}_p,d_q)(1+\delta\xi_{p,q}),
\enn
where $\delta$ is the noise level and $\xi_{p,q}$ is the standard normal distribution.
In the training stage, we choose $\delta=5\%$.

\subsubsection{Parameter setting for inversion algorithms}\label{S5_1_2}

For the parameters in \textit{Learned Projected Algorithm}, we choose the resolution
$N=64$, wave number $k=6$, $L=100$, stepsize $\mu=1$, $R=5$, $N_0=20$ and
the regularization parameters $\alpha_i:=5\times (0.2)^i,\;i=0,\dots,R-1$
(the choice of $\alpha_i$ follows the suggestion in \eqref{3}). In order to obtain
the well-trained neural network $\mathcal{P}_{\widehat{\Theta}}$
for \textit{Learned Projected Algorithm}, we train $\mathcal{P}_{\Theta}$ by minimizing
the error functions $\mathcal{E}_1$, $\mathcal{E}_2$, $\mathcal{E}_3$ and $\mathcal{E}$
with the epochs $t_1=100$, $t_2=20$, $t_3=20$ and $t_4=20$, respectively,
with using the Adam optimizer \cite{K14} with batch size $30$ and learning rate $10^{-3}$.
We also compare \textit{Learned Projected Algorithm} with \textit{Simplified Learned Projected Algorithm}
to illustrate the benefit brought by the training dataset $S_4$ (see Subsection \ref{S4_3}).
Here, the parameters in \textit{Simplified Learned Projected Algorithm} are the same as those
in \textit{Learned Projected Algorithm}, except that the learned projector
$\mathcal{P}_{\widehat{\Theta}}$ is replaced by $\mathcal{P}_{\Theta_3}$.

\subsubsection{Evaluation criterion for inversion algorithms}

In order to quantitatively evaluate the reconstruction performance of \textit{Learned Projected Algorithm},
we introduce an error function to measure the difference between the exact refractive index and
the approximate refractive index obtained by \textit{Learned Projected Algorithm}.
As mentioned before, for the exact contrast $m(x)$, $\boldsymbol{m}=(\boldsymbol{m}_{ij})\in\C^{N\times N}$
is the exact contrast matrix with $\boldsymbol{m}_{ij}=m(x_{ij})$ and
$\widehat{\boldsymbol{m}}=(\widehat{\boldsymbol{m}}_{ij})\in\C^{N\times N}$ is the output of
\textit{Learned Projected Algorithm} which is the approximation of $\boldsymbol{m}$.
Here, $x_{ij}$ ($i,j = 1,2,\ldots,N$) are the points introduced at the end of Section \ref{S2}.
Accordingly, $\boldsymbol{n}=(\boldsymbol{n}_{ij}):=\boldsymbol{m}+1$ is the discretization of
the refractive index $n(x)= m(x)+1$ and
$\widehat{\boldsymbol{n}}=(\widehat{\boldsymbol{n}}_{ij}):=\widehat{\boldsymbol{m}}+1$
is the approximation of $\boldsymbol{n}$.
Now we define the relative error function $R_e$ between $\boldsymbol{n}$ and $\widehat{\boldsymbol{n}}$
as follows
\begin{displaymath}
	R_e(\boldsymbol{n}, \widehat{\boldsymbol{n}}) := \sqrt{\dfrac{1}{N^2}
		\left[\sum_{i,j}\left|\dfrac{\boldsymbol{n}_{ij}
			- \widehat{\boldsymbol{n}}_{ij}}{\boldsymbol{n}_{ij}}\right|^2\right]}.
\end{displaymath}
Moreover, the reconstruction performance of other inversion algorithms carried out in
the numerical experiments will be evaluated in the same way as above.

\subsection{Performance of the proposed algorithm}\label{S5_2}

We train $\mathcal{P}_\Theta$ by using MNIST dataset \cite{D12}, which consists of $10$ handwritten
digits from $0$ to $9$. To be more specific, we randomly select $T=2000$ different digits from MNIST
dataset to represent the exact contrast matrices $\{\boldsymbol{m}^{(i)}\}_{i=1}^T$ of the unknown
inhomogeneous media such that $\|\boldsymbol{m}^{(i)}\|_{\max},\;i=1,\ldots,T$,
lie uniformly in the interval $[1,4]$. During the training process, $1800$ samples are used for
training $\mathcal{P}_\Theta$ with the training strategy in Section \ref{S4_2} and $200$ samples
are used to validate the training performance.

\begin{figure}[htpb]
	\centering
	\includegraphics[width=\textwidth]{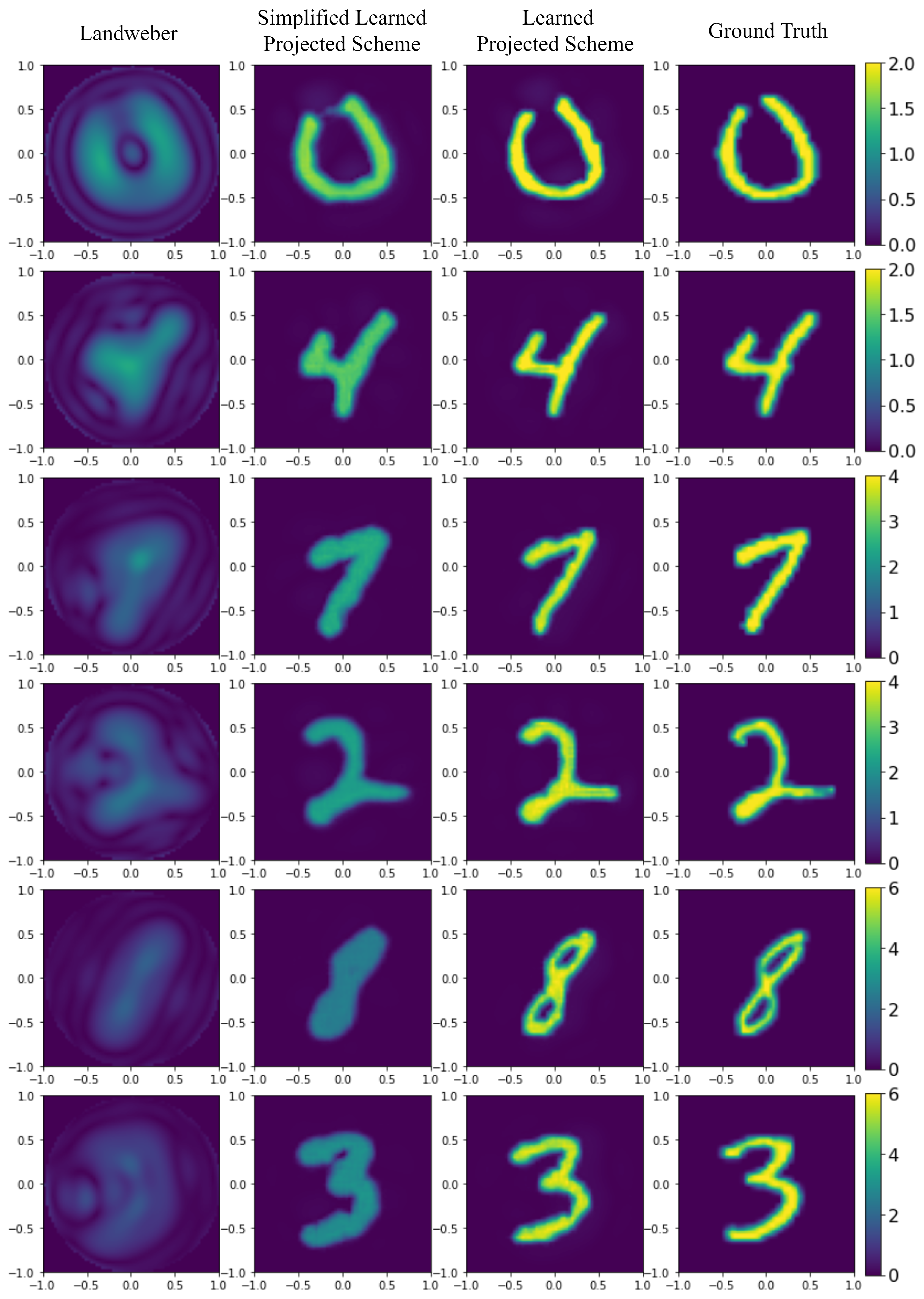}
	\caption{The initial guess generated by the Landweber method, reconstructions by
		\textit{Simplified Learned Projected Algorithm} and \textit{Learned Projected Algorithm}
		with exact contrast matrices generated from MNIST digits. Each row presents the reconstruction
		results and the ground truth for one sample.
	}\label{F3}
\end{figure}

First, in order to investigate the influence of values of different contrasts on
\textit{Learned Projected Algorithm}, we choose the exact contrast matrices generated from other MNIST
digits for evaluation and choose the noise level $\delta=5\%$, where the norm $\|\cdot\|_{\max}$ of each of these matrices is set to be $2,4$ or $6$.
We also carry out the same testing contrast matrices for \textit{Simplified Learned Projected Algorithm} with the noise level $\delta=5\%$.
The resulting reconstructed contrast matrices obtained by these two algorithms are illustrated in
Figure \ref{F3}. Each row of Figure \ref{F3} presents the reconstruction result of Landweber method,
the reconstruction result of \textit{Simplified Learned Projected Algorithm}, the reconstruction result
of \textit{Learned Projected Algorithm} and the ground truth for one sample.
In order to compare these two algorithms quantitatively, we consider three cases with the noise setting
$\delta=5\%$, which are denoted as Cases 1.1, 1.2 and 1.3.
In all three cases, we randomly generate $100$ samples from the MNIST dataset to represent the exact
contrast matrices.
For Cases 1.1, 1.2 and 1.3, we set the norm $\|\cdot\|_{\max}$ of each exact contrast matrix to be $2,4$
and $6,$ respectively. For these three cases, the second row and the third row in Table \ref{T1}
present the average values of the relative errors $R_e$ for the outputs of
\textit{Simplified Learned Projected Algorithm} and \textit{Learned Projected Algorithm}, respectively.
The reconstruction results in Figure \ref{F3} and Table \ref{T1} show that the proposed
\textit{Learned Projected Algorithm} could generate satisfactory results even for the case
when the values of the contrasts for testing are higher than those for training.
We believe this is because the feasible set $\mathcal{M}_R$ in \eqref{5} we learn for
\textit{Learned Projected Algorithm} mainly characterizes the a priori information of the shapes
of the unknown contrasts, since the training datasets for learned projector $\mathcal{P}_{\widehat{\Theta}}$
are obtained by using the normalization operator $\mathcal{N}$ (see Section \ref{S4_2}),
and hence the proposed \textit{Learned Projected Algorithm} could deal with unknown contrasts
with various values. Moreover, it can be seen in Figure \ref{F3} and Table \ref{T1} that our proposed
\textit{Learned Projected Algorithm} improves the initial guesses remarkably and outperforms
\textit{Simplified Learned Projected Algorithm}, which shows the advantages offered by
the deep learning method and the training dataset $S_4$.

\begin{figure}[htpb]
	\centering
	\includegraphics[width=\textwidth]{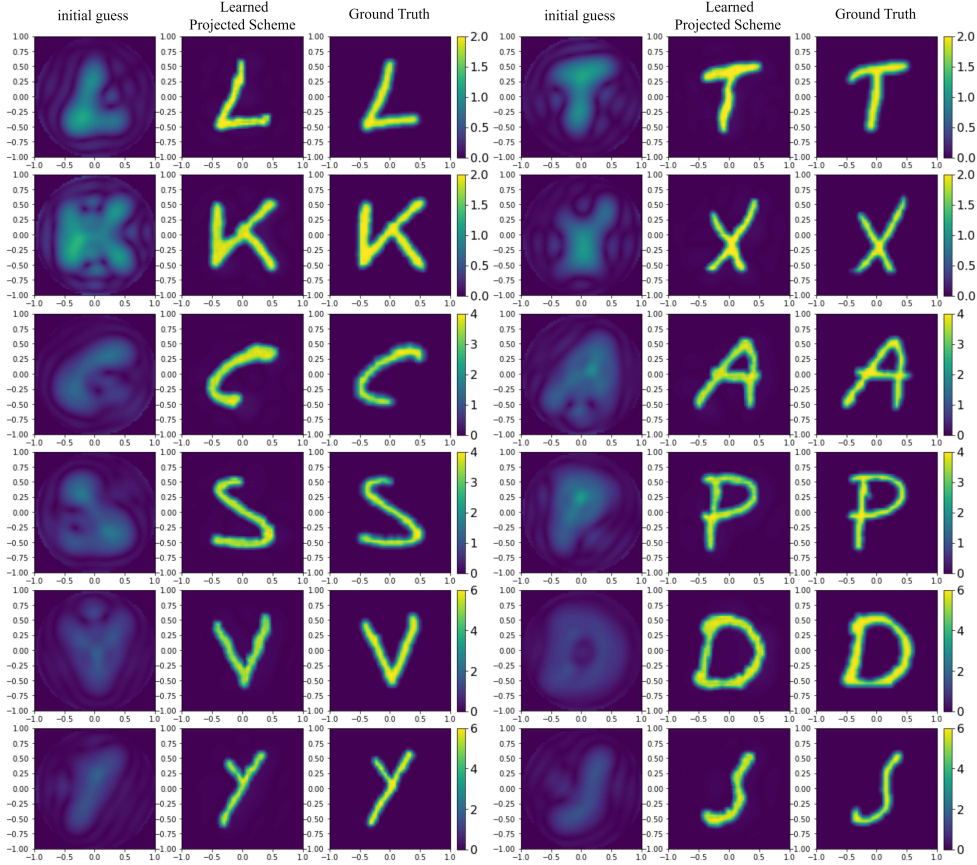}
	\caption{Reconstructions by \textit{Learned Projected Algorithm} with the exact contrast matrices
		generated from EMNIST dataset. Each row presents the reconstruction results and the ground truth
		for two samples.
	}\label{F4}
\end{figure}

Secondly, we evaluate the performance of the proposed algorithm on the EMNIST dataset \cite{C17}
by using exact contrast matrices generated from letters in the dataset, with a noise level of $\delta=5\%$.
We present the reconstruction results of \textit{Learned Projected Algorithm} in Figure \ref{F4},
where each row displays the reconstruction results and the ground truth for two samples.
Each sample in Figure \ref{F4} is presented with the reconstruction result of the Landweber method,
the reconstruction result \textit{Learned Projected Algorithm}, and the ground truth.
To quantitatively evaluate the performance of \textit{Learned Projected Algorithm} on EMNIST dataset,
we consider three cases with the noise setting $ \delta = 5\% $, which are denoted as Cases 2.1, 2.2 and 2.3.
In all three cases, we generate 100 samples from the MNIST dataset and 100 samples from the EMNIST dataset
to represent the exact contrast matrices. These 100 MNIST samples are the same as those used
in Cases 1.1, 1.2, and 1.3.
For Cases 2.1, 2.2 and 2.3, we set the norm $\|\cdot\|_{\max}$ of each exact contrast matrix to be 2, 4, and 6,
respectively. For these three cases, the second row and the third row in Table \ref{T2} present the
average values of the relative errors $ R_e $ for the outputs of \textit{Learned Projected Algorithm}
on MNIST dataset and EMNIST dataset, respectively.
It can be observed in Figure \ref{F4} and Table \ref{T2} that the proposed
\textit{Learned Projected Algorithm} has good performance on EMNIST dataset.
For this observation, it is reasonable to deduce that $\mathcal{M}_R$ in \eqref{5} we learn for
\textit{Learned Projected Algorithm} contains not only samples from MNIST dataset but also other
handwritten samples, which leads to satisfactory generalization ability of our algorithm.

\begin{table}[htbp]
\footnotesize
\caption{The average values of relative errors $ R_e $ for the outputs of
\textit{Simplified Learned Projected Algorithm} and \textit{Learned Projected Algorithm}
on the MNIST dataset.
}\label{T1}
\begin{center}
\begin{tabular}{cccc}
\hline
& Case 1.1 & Case 1.2 & Case 1.3\\
\hline
\textit{Simplified Learned Projected Algorithm} & 9.4\% & 21.3\% & 29.4\% \\
\hline
\textit{Learned Projected Algorithm} & 8.2\% & 17.9\% & 24.8\% \\
\hline
\end{tabular}
\end{center}
\end{table}

\begin{table}[htbp]
\footnotesize
\caption{The average values of relative errors $R_e$ for the outputs of \textit{Learned Projected Algorithm}
on the MNIST dataset and the EMNIST dataset.
}\label{T2}
\begin{center}
\begin{tabular}{cccc}
\hline
& Case 2.1 & Case 2.2 & Case 2.3\\
\hline
MNIST & 8.2\% & 17.9\% & 24.8\% \\
\hline
EMNIST & 10.1\% & 21.9\% & 28.7\% \\
\hline
\end{tabular}
\end{center}
\end{table}

\subsection{Sensitivity to noise}\label{S5_3}

\begin{figure}[h]
	\centering
	\includegraphics[width=\textwidth]{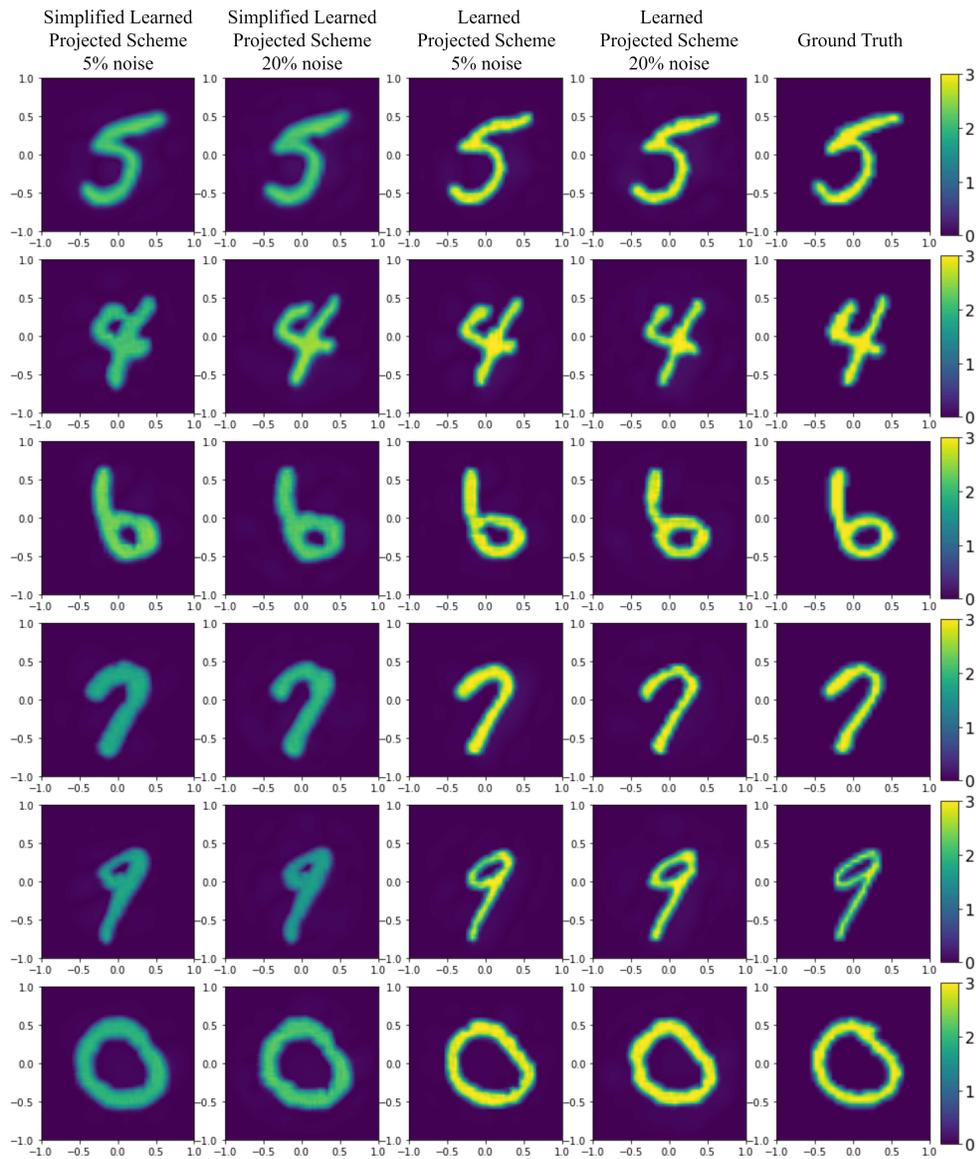}
	\caption{Reconstructions by \textit{Simplified Learned Projected Algorithm} and
		\textit{Learned Projected Algorithm} with different noise settings, where the exact contrast matrices
		are generated from MNIST dataset. Each row presents the reconstruction results and
		the ground truth for one sample.}\label{F5}
\end{figure}

In order to test the robustness of our algorithm, we test \textit{Simplified Learned Projected Algorithm}
and \textit{Learned Projected Algorithm} in the setting of different noises, where the parameters of
these two algorithms are the same as in Section \ref{S5_1_2}.
To do this, we consider two cases with different noise settings $\delta=5\%$ and $\delta=20\%$,
which are denoted as Cases 3.1 and 3.2, respectively.
In both two cases, we randomly generate $100$ samples from the MNIST dataset to represent the exact
contrast matrices with the norm $\|\cdot\|_{\max}$ to be $3$.
For these two cases, we present the average values of the relative errors $R_e$ for the outputs of
the above two algorithms in Table \ref{T3}.
Figure \ref{F5} presents the reconstruction results of several samples from these two cases.
Each row of Figure \ref{F5} presents the reconstruction results of the above two algorithms and
the ground truth for one sample.
The reconstruction results in Table \ref{T3} and Figure \ref{F5} demonstrate that the proposed
\textit{Learned Projected Algorithm} performs well in a high noise setting, which exactly meets
our expectation. As we mentioned before in Section \ref{S3}, a good choice of the feasible set
$\mathcal{M}_\mathcal{R}$ in \eqref{5} has the effect of denoising the noisy measurement.

\begin{table}[htbp]
\footnotesize
\caption{The average values of relative errors $ R_e $ for the outputs of
\textit{Simplified Learned Projected Algorithm} and \textit{Learned Projected Algorithm}
on the MNIST dataset, where the noise levels of Cases 3.1 and 3.2 are set to be $\delta= 5\%$
and $\delta = 20\%$, respectively.
}\label{T3}	
\begin{center}
\begin{tabular}{ccc}
\hline
& Case 3.1 & Case 3.2\\
\hline
\textit{Simplified Learned Projected Algorithm} & 15.5\% & 16.6\%\\
\hline
\textit{Learned Projected Algorithm} & 13.3\% & 15.1\%\\
\hline
\end{tabular}
\end{center}
\end{table}

\section{Conclusion}\label{S6}

In this paper, we considered the inverse problem of scattering of time-harmonic acoustic waves
from inhomogeneous media in two dimensions, including the high contrast case.
Such an inverse scattering problem is nonlinear and severely ill-posed, and thus certain
regularization method need to be used.
However, the selection of regularization functionals relies heavily on the unique demand of
individual application scenarios, which motivates us to learn certain a priori information of
the unknown scatterer from the ground truth data.
To this end, we reformulated the inverse scattering problem as a constrained minimization
problem \eqref{5} with an unknown constrained feasible region depending on the unknown
regularization parameter and functional.
We then proposed \textit{Learned Projected Algorithm} (see Section \ref{S4_3}) to solve \eqref{5},
which employs the learned projector $\mathcal{P}_{\widehat{\Theta}}$ to force the constraint
to be satisfied and the IRGNM algorithm to minimize the data-fitting term.
In our algorithm, the a priori information of the unknown scatterer is encoded in the learned
projector $\mathcal{P}_{\wi{\Theta}}$, which is a well-trained deep neural network that learns
the a priori information of the shape of the unknown scatterer directly from the ground truth data,
and is expected to provide high-quality initial guesses for IRGNM. Here, the a priori information
of the shape of the unknown scatterer has played a crucial role in dealing with the high contrast
case satisfactorily.

Various numerical experiments show that the \textit{Learned Projected Algorithm} performs
well for the inverse problem considered.
First, it is observed that the \textit{Learned Projected Algorithm} can generate satisfactory
reconstruction results for a wide range of contrasts with different values and even for the case
when the values of the contrasts for testing are slightly higher than those for training.
The reason for this may be due to fact that the learned projector acquires some a priori
information about the shapes of the unknown contrasts since the training datasets for the learned
projector are obtained by using the normalization operator $\mathcal{N}$ (see Subsection \ref{S4_2}).
Secondly, it is seen that \textit{Learned Projected Algorithm} has a good performance on the EMNIST
dataset though we only train $\mathcal{P}_{{\Theta}}$ on the MNIST dataset.
This suggests that the feasible region we learned for \textit{Learned Projected Algorithm} contains
not only the samples from the MNIST dataset but also other handwritten samples, leading to
satisfactory generalization ability for our algorithm.
Thirdly, as mentioned in Section \ref{S3}, a good choice of the feasible region has the effect
of denoising the noise measurement, which is consistent with our experiments that the performance
of \textit{Learned Projected Algorithm} does not degrade significantly as the noise level increases.
However, it is observed in the numerical experiments that the reconstruction results of
\textit{Learned Projected Algorithm} become worse when the value of the exact contrast is large.
%This may be due to two reasons. The first reason is
One of the reasons may be due to the fact that for the case when the value of the exact contrast is
large the numerical solution of the corresponding scattering problem in each iterative step is not
accurate enough, which will deteriorate the final reconstruction results.
%The second reason is that the deep learning method in our algorithm is lack of interpretability
%(that is, the underlying a priori information acquired by the learned projector is not well-understood)
%and thus it is hard for us to improve the proposed algorithm.
%Therefore, certain improvements still need to be investigated.
%Moreover, it is interesting
It is also interesting to extend our method to the case of seismic imaging, which will be considered
as a future work.

\section*{Acknowledgments}

This work was partly supported by the National Key R\&D Program of China (2018YFA0702502),
Beijing Natural Science Foundation (Z210001), the NNSF of China (12271515) and Youth Innovation
Promotion Association of CAS.

\providecommand{\href}[2]{#2}
\providecommand{\arxiv}[1]{\href{http://arxiv.org/abs/#1}{arXiv:#1}}
\providecommand{\url}[1]{\texttt{#1}}

\end{document}